\documentclass[11pt]{article}

\usepackage{amsmath, amssymb, amscd}

\def\eqref#1{(\ref{#1})}

\newcommand{\arrow}{{\:\longrightarrow\:}}
\newcommand{\Z}{{\Bbb Z}}

\newcommand{\R}{{\Bbb R}}

\def\1{\sqrt{-1}\:}
\newcommand{\restrict}[1]{{\left|_{{\phantom{|}\!\!}_{#1}}\right.}}

\renewcommand{\phi}{\varphi}
\renewcommand{\epsilon}{\varepsilon}
\renewcommand{\geq}{\geqslant}

\newcommand{\Aut}{\operatorname{Aut}}

\newcommand{\comment}[1]{{}}

\def\blacksquare{\hbox{\vrule width 4pt height 4pt depth 0pt}}
\def\endproof{\blacksquare}

\makeatletter

\@ifundefined{Bbb}
     {\newcommand{\Bbb}[1]{{\mathbb #1}}}%
{}

\newcommand{\ps@verbit}{%
  \renewcommand{\@oddhead}{%
          \scriptsize
          {Structure theorem for Vaisman manifolds}
          \hfil\tiny {L. Ornea, M. Verbitsky, 18.05.2003}}
  \renewcommand{\@evenhead}{\@oddhead}
  \renewcommand{\@oddfoot}{\hfil\thepage\hfil}
  \renewcommand{\@evenfoot}{\@oddfoot}}
 
\newcounter{Mycounter}[section]
\newcounter{lemma}[section]
\renewcommand{\thelemma}{\noindent{Lemma \thesection.\arabic{lemma}}}
\newcommand{\lemma}{%
     \setcounter{lemma}{\value{Mycounter}}
     \refstepcounter{lemma}
     \stepcounter{Mycounter}
     {\bf \thelemma:\ }}

\newcounter{claim}[section]
\newcounter{sublemma}[section]
\newcounter{corollary}[section]
\newcounter{theorem}[section]
\renewcommand{\thetheorem}{\noindent{Theorem \thesection.\arabic{theorem}}}
\newcommand{\theorem}{%
     \setcounter{theorem}{\value{Mycounter}}
     \refstepcounter{theorem}
     \stepcounter{Mycounter}
     {\bf \thetheorem:\ }}

\newcounter{theoremm}
\renewcommand{\thetheoremm}{\noindent{Structure Theorem 
}}
\newcommand{\theoremm}{%
     \setcounter{theorem}{\value{Mycounter}}
     \refstepcounter{theorem}
     \stepcounter{Mycounter}
     {\bf \thetheoremm:\ }}

\newcounter{conjecture}[section]
\newcounter{proposition}[section]
\renewcommand{\theproposition}
       {\noindent{Proposition \thesection.\arabic{proposition}}}
\newcommand{\proposition}{%
     \setcounter{proposition}{\value{Mycounter}}
     \refstepcounter{proposition}
     \stepcounter{Mycounter}
     {\bf \theproposition:\ }}

\newcounter{definition}[section]
\renewcommand{\thedefinition}
       {\noindent{Definition~\thesection.\arabic{definition}}}
\newcommand{\definition}{%
     \setcounter{definition}{\value{Mycounter}}
     \refstepcounter{definition}
     \stepcounter{Mycounter}
     {\bf \thedefinition:\ }}

\newcounter{example}[section]
\newcounter{remark}[section]
\renewcommand{\theremark}{\noindent{Remark \thesection.\arabic{remark}}}
\newcommand{\remark}{%
     \setcounter{remark}{\value{Mycounter}}
     \refstepcounter{remark}
     \stepcounter{Mycounter}
     {\bf \theremark:\ }}

\newcounter{problem}[section]
\newcounter{question}[section]
\@addtoreset{equation}{section}
\@addtoreset{footnote}{section}
\makeatother

\begin{document}

\begin{center}
{\LARGE\bf
Structure theorem for compact Vaisman manifolds 
}
\\[4mm]
Liviu Ornea\footnote{Liviu Ornea is a member of EDGE, Research
Training Network HRPN-CT-2000-00101, supported by the European Human
Potential Programme.}, Misha Verbitsky,\footnote{Misha Verbitsky is 
an EPSRC advanced fellow 
supported by CRDF grant RM1-2354-MO02 and EPSRC grant  GR/R77773/01.

{\it Keywords and phrases~:} Locally conformal K{\"a}hler manifold,
Sasakian manifold, parallel
Lee form, Gauduchon metric, weight bundle, monodromy.

2000 {\it MSC}~: 53C55, 53C25.}\\[4mm]

\end{center}

{\small 
\hspace{0.15\linewidth}
\begin{minipage}[t]{0.7\linewidth}
{\bf Abstract} \\
A locally conformally K{\"a}hler (l.c.K.) manifold is a complex manifold
admitting a K{\"a}hler covering $\tilde M$, with each deck transformation 
acting by K{\"a}hler homotheties.
A compact l.c.K. manifold is Vaisman if it admits a holomorphic
flow acting by non-trivial homotheties on $\tilde M$.
We prove a structure theorem for compact Vaisman manifolds.
Every compact Vaisman manifold $M$ is fibered over a circle,
the fibers are Sasakian, the fibration is locally trivial,
and $M$ is reconstructed from the Sasakian structure 
on the fibers and the monodromy automorphism induced
by this fibration. This construction is canonical and 
functorial in both directions.
\end{minipage}
}

\section{Introduction}

Let  $(M,J,g)$  be a connected complex 
manifold of complex dimension
$n \geq 2$, equipped with a Hermitian metric.  
We shall denote by $\omega$ its fundamental 
two-form given by $\omega(X,Y)=g(X,JY)$.

\hfill

\definition
$(M,J,g)$ is called {\em
 locally conformal   K{\"a}hler},  {\em l.c.K.} for short,  if there
 exists an open cover $\mathcal{U}=\{U_{\alpha}\}$ such that each
 locally defined metric 
$g_{\alpha}=e^{-f_\alpha}g_{\restrict{U_\alpha}}$ 
is K{\"a}hler for some smooth  
      function
       $f_\alpha$ on $U_\alpha$. 
\hfill

Equivalently,   $(M,J,g)$ is l.c.K. if and only if there exists a
{\em closed} one-form $\theta$ such that
\begin{equation}\label{1}
d\omega=\theta\wedge  \omega.
\end{equation}

\noindent Locally, $\omega_{|U}=df_U$. Note also that on
complex surfaces, the equation \eqref{1} implies $d\theta =0$. 

The globally defined one-form $\theta$ is called {\em the Lee form} and its metrically
equivalent (with respect to $g$) vector field $V=\theta^\sharp$ 
is called {\em the Lee
vector field}.

An equivalent definition of  l.c.K. manifolds is obtained
\emph{via} the universal Riemannian cover. Let $p~:\tilde M\rightarrow
M$ be the covering map and denote also by $J$ the lifted complex structure. The lifted
metric $\tilde g$ is globally conformal with a K{\"a}hler metric $g_k$, because $p^*\theta$
is exact. The fundamental group of $M$ will then act by holomorphic
conformal transformations with respect to $(J,\tilde h)$. But
conformal transformations of a symplectic form (in real dimension at
least $4$) are in fact homotheties (note that this remains true for the deck group of any covering). The converse is also true,
namely: $(M,g,J)$ is l.c.K. if and only if the universal Riemannian
covering admits a K{\"a}hler form $\omega_k$ with respect to which $\pi_1(M)$
acts by holomorphic homotheties (cf. \cite{_Dragomir_Ornea_}). 

For later use, let $\mathcal{H}(\tilde M, \tilde h,J)$ be the above group of
holomorphic homotheties and let 
\begin{equation}\label{_homotheties_Equation_}
\rho~:\mathcal{H}(\tilde M, \tilde h,J)\arrow \R^{>0}
\end{equation}
be the homomorphism associating to each homothety
its scale factor.

It will be useful in the sequel to regard l.c.K. geometry in the
framework of Weyl geometry (cf. \emph{e.g.} \cite{_Gauduchon_Crelle_}). Let $L\rightarrow M$ be the weight bundle,
\emph{i.e.} the bundle of scalars of weight
1 associated to the bundle of linear frames of $(M,g)$ by the
representation $GL(n,\R)\ni A\mapsto |\det A|^{\frac 1n}$. The Lee
form $\theta$ can be viewed as a connection form in $L$. The induced
Weyl connection is clearly flat. We shall refer to the holonomy of
this connection as to the monodromy of $L$.  

A strictly smaller class of l.c.K. manifolds is the one formed by those with
parallel (with respect to the Levi Civita connection) Lee form. We
call  them {\em Vaisman manifolds}, as I. Vaisman was the first to
study them sytematically (under the name of generalized Hopf manifolds
 \cite{_Vaisman:Torino_}, \cite{_Vaisman:Dedicata_}, a name which later proved to be inappropriate). 
On such a  manifold, the length of
the Lee vector field is constant and we shall always assume it is
nonzero. Hence, in what follows, we shall normalize and consider that on a
 Vaisman
manifold $|v|=1$. 

Note that on a compact l.c.K. manifold, the metric with parallel Lee form, if
it exists, is unique up to homothety in its conformal class and
coincide with the Gauduchon (standard) metric, \cite{mpps}. Hence, on
a compact Vaisman manifold, we may always assume that we work with a
metric with harmonic Lee form of length $1$. Moreover,
\emph{loc. cit.}, the group $\mathop\mathrm{Aut}(M)$ of l.c.K. automorphisms of a compact
Vaisman manifold coincides with the isometry group of the Gauduchon
metric, thus being a compact Lie group.

To describe examples of Vaisman manifolds we need to recall the notion
of Sasakian manifold (see \cite{blair} and  \cite{_Boyer_Galicki_} for
a survey and references on Sasakian geometry). We need first the following

\hfill

\definition
Let $(X,g)$ be a Riemannian manifold.
A {\bf cone} over $X$ is a Riemannian manifold
$C(X):= (X \times \R^{>0}, dt^2 + t^2 g )$,
where $t$ is the parameter on $\R^{>0}$. 
For any $\lambda\in \R^{>0}$,
the map 
\[ \tau_\lambda:\; C(X) \arrow C(X), \ \  (x, t) \arrow (x, \lambda t)
\]
multiplies the metric by $\lambda^2$.

\hfill

\definition
Let $X$ be a Riemannian manifold. A Sasakian  
structure on $X$ is a complex structure on $C(X)$ satisfying the following
\begin{description}
\item[(i)] The metric on $C(X)$ is K{\"a}hler.
\item[(ii)] The map $\tau_\lambda:\; C(X) \arrow C(X)$ is holomorphic,
for all $\lambda\in \R^{>0}$.
\end{description}

\hfill

Let $X$ be a Sasakian manifold, and $C(X)$ its K{\"a}hler cone. Given a number
 $q\in \R$, $q>1$, consider an equivalence relation
$\sim_q$ on $C(X) = X\times \R^{>0}$ generated
by $(x,t) \sim (x, qt)$. Since the
map $(x,t) \arrow (x, qt)$ mupltiplies the metric
by $q^2$, the quotient $M = C(X)/\sim_q$
is an l.c.K. manifold. Moreover,
$M$ is a Vaisman manifold, with the Gauduchon
metric provided by an isomorphism
$M \cong X \times S^1$.

Conversely, it follows from \cite{_Kamishima_Ornea_} (see also 
\cite{_Gini_Ornea_Parton_}, \cite{_Verbitsky:LCHK_}) that 
 any Vaisman manifold $(M,g,J)$ with exact Lee form 
is a cone of a naturally defined 1-Sasakian manifold. Indeed, if
$\theta = df$,  since $|\theta| = const$,
$f:\; M \arrow \R$ is a Riemannian submersion. The Lee flow
$\psi_t$ is the gradient flow of $f$. Therefore, $\psi_t$
is compatible with $f$ in the following sense:
\[ f (\psi_t x) = f (x) +t .
\]
As the Lee field is Killing, $\psi_t$ is formed by isometries.
Therefore, $\psi_t$ induces (locally in $M$) 
a trivialization of the fibration $f:\; M \arrow \R$
identifying its fibers. 

Let $X$ be the fiber of $f$, and let $\rho:\; M \arrow X$ 
the corresponding trivialization map. The map
$\rho$ is well defined everywhere in $M$ if
$\psi_t$ is defined for all $t$.

The manifold $M$ is now isomorphic to $X \times \R$.
Using exponential map, we identify $X \times \R$ with
$X \times \R^{>0}$. The cone metric on $M$ is written 
as $e^{-f} g$ and it is K{\"a}hler. By definition,
this implies that $X$ is a Sasakian manifold.

More generally, 
let $\phi:\; X \arrow X$ be an automorphism of
a Sasakian structure. The map $(x, t) \stackrel{\phi_q} \arrow (\phi(x), qt)$
is compatible with the complex structure and multiplies
the metric by $q^2$. Therefore, the quotient
$M_{\phi, q}$ of $C(X)$ by the
corresponding equivalence relation
$\sim_{\phi,q}$ is called the suspension of $\phi$ over the circle of
length $2\pi q$ and is a l.c.K. manifold.
The following proposition is implied by \cite{_Kamishima_Ornea_}
(see also \cite{_Gini_Ornea_Parton_}). 

\hfill

\proposition
Assume that the Sasakian manifold $X$ 
is compact, and $M_{\phi, q}:= C(X)/\sim_{\phi,q}$ an l.c.K. manifold
constructed as above. Then 
$M_{\phi, q}$ is a Vaisman manifold.

\hfill

The most common examples of this kind are the Hopf surfaces with the
metrics constructed in \cite{go}, \cite{_Belgun_}, and generalized in
\cite{_Kamishima_Ornea_}.

The aim of this paper is to prove the following 
converse to the above proposition, which is proven
for Einstein-Weyl locally conformally K{\"a}hler manifolds
in \cite{_Verbitsky:LCHK_}:

\hfill
 
\theoremm \label{_structure_Theorem_}
Let $M$ be a compact Vaisman manifold. Then $M$
admits a canonical Riemannian submersion $p:\; M \arrow S^1$ to a circle,
and the fibers of $p$ are isomorphic as Riemannian manifolds. 
Moreover, for all $s \in S$, the manifold $X=p^{-1}(s)$
is equipped with a natural 1-Sasakian structure
and a 1-Sasakian sutomorphism $\phi$, such that
$M$ is isomorphic to the Vaisman manifold
$M_{\phi, q}:= C(X)/\sim_{\phi,q}$ constructed 
above.

\section{Proof of the structure theorem}

Let $(M,g,J)$ be a compact Vaisman manifold such that $g$ is the
Gauduchon metric. The associated Lee form $\theta$ being harmonic, $M$
is foliated by minimal hypersurfaces whose normal trajectories are
geodesics. 

To obtain a submersion over $S^1$ we proceed as follows. Fix a
point $x_0\in M$. Consider a smooth path $\gamma:\; [0,1] \arrow M$,
$\gamma(0)=x_0$, $\gamma(1)=x$. Denote by $R(\gamma, x)$ the integral
\[ R(\gamma,x) := \int_\gamma \theta.
\]
Since $\theta$ is closed, $R(\gamma,x )$ depends only from the
homotopy class of $\gamma$. 

Let $\gamma_1, ... , \gamma_k$ be the 
generators of $H^1(M, \Z)$. Denote by $\alpha_1, ... \alpha_k$ the
periods of $\theta$,
\[ \alpha_i := \int_{\gamma_i} \theta. \]
Clearly, $R(\gamma,x)$ defines a function
\[ p:\; M \arrow \R / \langle \alpha_1, ... \alpha_k\rangle, \]
where $\langle \alpha_1, ... \alpha_k\rangle\subset \R$ is the
abelian group generated by $\alpha_1, ... \alpha_k$ 
and 
$$ \R / \langle \alpha_1, ... \alpha_k\rangle$$ 
is the quotient
group. 

If
$\langle \alpha_1, ... \alpha_k\rangle\cong \Z$,
then we have the desired submersion over $S^1$. The above condition is
certainly satisfied if the monodromy of the flat connection associated
to $\theta$ in $L$ is $\Z$. 
 Indeed, let $\gamma$ be a closed loop in $M$. By definition of
$L$, the monodromy of $L$ along $\gamma$ is equal to 
$\int_{\gamma} \theta$. If the monodromy of $L$ is isomorphic
to $\Z$, then the periods
\[ \alpha_i := \int_{\gamma_i} \theta \]
are proportional (with integer coefficient)
to $\alpha_0\in \R$. Then the function 
$p:\; M \arrow \R / \langle \alpha_1, ... \alpha_k\rangle$
maps $M$ to $\R / \Z \alpha_0 =S^1$.
The differential of $p$ is equal $\theta$.
Passing to a covering, we obtain a map
$f:\; \tilde M \arrow \R$, $df =\theta$,
$\tilde M / Z \cong M$. As we have shown above,
$\tilde M$ is a cone over a Sasakian manifold:
$\tilde M \cong C(X)$. The monodromy map
$\mu:\; \tilde M \arrow \tilde M$ preserves the complex
structure and the Lee field and multiplies the K{\"a}hler
form by a constant. Therefore, $\gamma$ is induced
by a Sasakian automorphism $\phi:\; X \arrow X$.
We obtain that $M \cong \tilde M /\mu \cong C(X)/\sim{q, \phi}$.

Hence, the structure
theorem will be implied by the following:

\hfill

\theorem\label{_weight_mono_Z_Corollary_}
Let $M$ be a compact Vaisman manifold, $L$ its weight bundle
and $\Gamma \subset \R^{>0}$ its monodromy group. Then 
$\Gamma \cong \Z$.

\hfill

For the proof,  we need to look at the group generated by the Lee flow and
at its lift on $\tilde M$. Let $\pi:\; \tilde M\arrow M$ 
the covering associated with the monodromy group of $L$
and  denote by $\Aut(\tilde M, M)$ 
the group of all conformal automorphisms of $\tilde M$
which make the following diagram commutative: 

\begin{equation}
\begin{CD}
 \tilde M@>{\tilde f}>> \tilde M \\
@V{\pi}VV  @VV{\pi}V              \\
M@>{f}>>  M 
\end{CD}
\end{equation}

$\Aut(\tilde M, M)$ is nonempty, as it contains the Lee flow $\tilde \Psi_t$. 

Let 
\begin{equation}\label{_M_automo_forgetful_Equation_} 
 \Phi:\; \Aut (\tilde M, M)\arrow \Aut (M)
\end{equation}
be the forgetful map. 
Recall that in \eqref{_homotheties_Equation_}
we defined the homomorphism $\rho$ mapping the group
of holomorphic homotheties of $\tilde M$ to $\R^{>0}$.
Restricting $\rho$ to $\Aut(\tilde M, M)$, we obtain
a map $\Aut(\tilde M, M) \arrow \R^{>0}$.

\hfill

\lemma \label{_aut(tilde_M,M)_inje_Lemma_}
In the above assumptions, the product map 
\[ 
   \Aut (\tilde M, M)\stackrel {\Phi \times \rho} \arrow \Aut (M)\times\R^{>0}
\]
is injective.

\hfill

\noindent{\bf Proof:} 
Let $\tilde f:\; \tilde M \arrow \tilde M$ be an automorphism
such that $\Phi(\tilde f) = id$. Then $\tilde f$ is a deck transform
of $\tilde M$. However, the group $\Gamma$ of the deck transforms
of $\tilde M$ is identified with the monodromy of the bundle $L$
of K{\"a}hler forms on $M$. Therefore, the natural map
\[ 
  \rho:\; \Gamma \arrow \R^{>0}
\]
is injective.
\endproof

Now let $G$ be the smallest Lie subgroup of $\Aut(M)$
containing the Lee flow $\psi_t$. By results in 
\cite{_Kamishima_Ornea_}, $G$ is a connected, compact abelian
group of holomorphic isometries of $(M,g)$. Let $\tilde G=\Phi^{-1}(G)$.
The main step in the proof of the structure theorem is the following:
  
\hfill

\theorem\label{_tilde_G_structure_Theorem_}
In the above assumptions, we have a (non-canonical) isomorphism
\[ \tilde G \cong \R \times (S^1)^k,
\]
where $k = \dim G -1$.

\hfill

\noindent{\bf Proof:}
Since $G$ is a compact abelian Lie group,
we have an isomorphism $G \cong (S^1)^{k+1}$.
By construction, $\tilde G$ is a covering of $G$.
In \ref{_aut(tilde_M,M)_inje_Lemma_}, we constructed
a natural monomorphism $\tilde G \arrow G \times \R$.
Denote by $\tilde G_0$ the connected component of 
$\tilde G$. Since $G$ is compact, the composition
\begin{equation} \label{_tilde_G_to_R_Equation_}
  \tilde G \arrow G \times \R \arrow \R
\end{equation}
is proper. Therefore, $\tilde G_0$ is isomorphic
to $\R \times (S^1)^k$ or to $(S^1)^{k+1}$.
The Lee flow $\tilde \phi_t$ multiplies the K{\"a}hler
form $\omega_k$ by a number, hence
the composition \eqref{_tilde_G_to_R_Equation_} 
is surjective. This implies that $\tilde G_0$
is non-compact. We obtain an isomorphism 
$\tilde G_0\cong \R \times (S^1)^k$.

Denote by $H$ the group $\tilde G / \tilde G_0$ of
connected components of $\tilde G$.
We first show that $H$ is finite, as follows.

Consider the restriction  $\rho:\; \tilde G \arrow \R^{>0}$,
where $\rho$ is the homomorphism defined in \eqref{_homotheties_Equation_}. 
Let 
\[ \tilde G_{\rho=1}:= \rho^{-1}(1)\]
be the corresponding subgroup (formed by K{\"a}hler isometries with
respect to $\omega_k$) of $\tilde G$. 
Clearly, $\rho$ is obtained by exponentiating
\eqref{_tilde_G_to_R_Equation_}.
Since $\rho:\; \tilde G \arrow \R$ is proper, 
$\tilde G_{\rho=1}$ is compact. To prove that $H$ is finite,
we only need to show  that every connected component
of $\tilde G$ meets $\tilde G_{\rho=1}$.

Let $\tilde h \in \tilde G$ be an arbitrary element,
$\rho(\tilde h) = e^t$. By construction,
$\rho(\tilde \phi_{-t}) = e^{-t}$, where
$\tilde \phi_t \in \tilde G_0$ is the
Lee flow of $\tilde M$. Then 
$\rho(\tilde h \cdot \tilde \phi_{-t})=1,$ and
\[ \tilde h \cdot \tilde \phi_{-t} \in \tilde G_{\rho=1}.
\]
On the other hand, $\tilde h$ and 
$\tilde h \cdot \tilde \phi_{-t}$ belong to the same
connected component of $\tilde G$. We obtained that
every connected component of $\tilde G$ meets
$\tilde G_{\rho=1}$, hence the group $H$
of connected components of $\tilde G$ is compact.

To finish the proof of \ref{_tilde_G_structure_Theorem_},
it remains to show that the group $H$ is trivial.
By the structure theorem for abelian Lie groups,
we have
\[ \tilde G \cong \R \times (S^1)^k\times H.
\]
This gives an embedding $H \hookrightarrow \tilde G$.
Denote by $\tilde H\subset \tilde G$ its image.
By definition, $\tilde H$ is a kernel
of a projection map $\tilde G \arrow \tilde G_0$. 
Since $\tilde G_0$ is mapped surjectively onto $G$,
the group $\tilde H$ is mapped to zero under
the forgetful map $\tilde G \arrow G$. Therefore,
$\tilde H$ is a subgroup of the monodromy group 
$\Gamma$. Since $\Gamma$ is the monodromy of the
real vector bundle $L$, we have $\Gamma\subset \R$,
and $\Gamma$ has no torsion subgroups. This
implies that $\tilde H=0$. We proved \ref{_tilde_G_structure_Theorem_}.
\endproof

\hfill

Now the proof of the structure theorem can be completed as follows. By definition of $\tilde G$,
the deck transformation group $\Gamma\subset \Aut(\tilde M, M)$ 
belongs to $\tilde G$. 
Note that $\Gamma = \ker \Phi$. Then, since    
$\tilde G \cong \R \times (S^1)^k$ and
$G \cong (S^1)^{k+1}$, $\ker \Phi \cong \Z$.

\hfill

\remark The proof used the  hypothesis of $(M,g,J)$ being a Vaisman manifold
under the equivalent statement that the Lee field has constant length
and its flow consists of holomorphic isometries of the Gauduchon metric. The example of the
Tricerri metric on the Inoue surface $S_m$ (see \cite{_Belgun_},
\cite{_Dragomir_Ornea_}) shows that harmonicity and constant length of
the Lee form are not enough, even when $b_1=1$. Indeed, Tricerri's metric has a harmonic
Lee form of (non-zero) constant length, so the fibration over $S^1$
exists, but the fibres are not Sasakian.

\hfill

\noindent{\bf Acknowledgements~:} This work was finished during a
visit to Ecole Polytechnique (Palaiseau). The authors thank the
colleagues in the Centre de Math{\'e}matiques for warm hospitality.

{\small

\hfill

\noindent {\sc Liviu Ornea\\
University of Bucharest, Faculty of Mathematics, 14
Academiei str., 70109 Bucharest, Romania.}\\
\tt lornea@imar.ro

\hfill

\noindent {\sc Misha Verbitsky\\
University of Glasgow, Department of Mathematics, 15
  University Gardens, Glasgow, United Kingdom.}\\
\tt verbit@maths.gla.ac.uk, \ \  verbit@mccme.ru 
}


\begin{thebibliography}{XXX}

\bibitem[Bel]{_Belgun_}
F.A. Belgun,
{\em On the metric structure of non-K{\"a}hler 
complex surfaces},
Math. Ann. {\bf 317} (2000), no. 1, 1--40.

\bibitem[Bes]{_Besse:Einst_Manifo_} 
A. Besse, 
 {\em Einstein Manifolds}, Springer-Verlag, New York, 1987.

\bibitem[Bl]{blair} D.E.  Blair, Riemannian geometry of contact and
symplectic manifolds, Progress in Math. {\bf 203}, Birkh{\"a}user,
Boston, Basel, 2002.

\bibitem[BG1]{_Boyer_Galicki_} 
C.P. Boyer, 
K. Galicki,  {\em 3-Sasakian Manifolds},
hep-th/9810250, also published in
Surveys Diff.Geom. {\bf 7} (1999) 123--184.

\bibitem[DO]{_Dragomir_Ornea_}
S. Dragomir,  L. Ornea,  
{\em Locally conformal K{\"a}hler geometry},
Progress in Mathematics, 155. 
Birkh{\"a}user, Boston, MA, 1998. 

\bibitem[G1]{_Gauduchon_1984_} 
P Gauduchon, 
 {\em La 1-forme de torsion d'une variete 
hermitienne compacte}, Math. Ann., {\bf 267} (1984), 495--518.

\bibitem[G2]{_Gauduchon_Crelle_} P. Gauduchon, \emph{Structures de
    Weyl-Einstein, vari{\'e}t{\'e}s de twisteurs et vari{\'e}t{\'e}s de type
    $S^1\times S^3$}, J. Reine Angew. Math., {\bf 469} (1995), 1--50.

\bibitem[GO]{go} P. Gauduchon, L. Ornea,
{\em Locally conformally K{\"a}hler metrics on Hopf surfaces},
Ann. Inst. Fourier  {\bf 48} (1998), 1107-1127.

\bibitem[GOP]{_Gini_Ornea_Parton_}
R. Gini, L.  Ornea, M. Parton, 
{\em Locally conformal K{\"a}hler reduction}, math.DG/0208208,
25 pages.


\bibitem[KO]{_Kamishima_Ornea_} 
Y. Kamishima, L. Ornea,
{\em Geometric flow on compact locally conformally Kahler manifolds},
math.DG/0105040, 21 pages.

\bibitem[MPPS]{mpps} A.B. Madsen, H. Pedersen, Y.S. Poon and A. Swann,
  \emph{Compact Einstein-Weyl manifolds with large symmetry group},
  Duke Math. J. {\bf 88} (1997), 407--434.

\bibitem[Sa]{_Sasaki_} Sasaki,  S.,
{\em 
On differentiable manifolds with certain 
structures which are closely related to
almost contact structure}, 
Tohoku Math. J. {\bf 2} (1960), 459--476.


\bibitem[V1]{_Vaisman:Dedicata_} 
I. Vaisman, {\em  Generalized Hopf 
manifolds}, Geom. Dedicata {\bf 13} (1982), no. 3, 231--255. 

\bibitem[V2]{_Vaisman:Torino_} 
I, Vaisman,  
{\em A survey of generalized Hopf manifolds}, 
Rend. Sem. Mat. Torino, Special issue (1984), 205--221.


\bibitem[Ve]{_Verbitsky:LCHK_}
M. Verbitsky,
{\em  Vanishing theorems for locally 
conformal hyperk{\"a}hler manifolds}, 2003, 
math.DG/0302219, 41 pages .


\end{thebibliography}
\end{document}